\documentclass[11pt,a4paper]{article}

\usepackage[T1]{fontenc}
\usepackage[utf8]{inputenc}
\usepackage{lmodern}

\usepackage[left=0.7cm, right=0.7cm, top=1.4cm, bottom=1.4cm]{geometry}

\usepackage{amsmath,amssymb,amsthm,mathtools,bm,mathrsfs,stmaryrd}
\numberwithin{equation}{section}

\usepackage{enumitem}
\setlist{itemsep=2pt,topsep=3pt,parsep=0pt}
\usepackage[dvipsnames]{xcolor}
\definecolor{oiblue}{HTML}{0072B2}
\definecolor{oiorange}{HTML}{E69F00}
\definecolor{oigreen}{HTML}{009E73}
\definecolor{oivermilion}{HTML}{D55E00}
\definecolor{oisky}{HTML}{56B4E9}
\definecolor{oipurple}{HTML}{CC79A7}
\definecolor{oiyellow}{HTML}{F0E442}

\usepackage{booktabs,array,tabularx,multirow,adjustbox}
\newcolumntype{Y}{>{\centering\arraybackslash}X}   

\usepackage{graphicx}
\graphicspath{{./}{figures/}}
\usepackage{float,placeins}
\usepackage[font=small,labelfont=bf,skip=4pt]{caption}
\usepackage{subcaption}
\usepackage{tikz}
\usetikzlibrary{arrows.meta,positioning,calc,patterns,decorations.pathreplacing,shapes.geometric}
\usepackage{pgfplots}
\pgfplotsset{compat=newest}
\pgfplotscreateplotcyclelist{mathpaths}{
	{oiblue,mark=*,thick},{oivermilion,mark=square*,thick},{oigreen,mark=triangle*,thick},
	{oiorange,mark=diamond*,thick},{oipurple,mark=o,thick},{oisky,mark=x,thick}}
\pgfplotsset{
	width=0.62\textwidth, height=0.36\textwidth,
	cycle list name=mathpaths,
	every axis/.append style={grid=major, grid style={gray!25}, tick label style={font=\small},
		label style={font=\small}, legend style={font=\small, draw=none, fill=none}}}
\newcommand{\fig}[3]{\begin{figure}[tbp]\centering\includegraphics[width=#2]{#1}\caption{#3}\end{figure}}
\newcommand{\optfig}[3]{\IfFileExists{#1}{\fig{#1}{#2}{#3}}{}}

\usepackage{algorithm}
\usepackage[noend]{algpseudocode}
\usepackage{listings}
\newtheorem{theorem}{Theorem}[section]
\newtheorem{proposition}[theorem]{Proposition}
\newtheorem{lemma}[theorem]{Lemma}
\newtheorem{corollary}[theorem]{Corollary}
\newtheorem{conjecture}[theorem]{Conjecture}
\theoremstyle{definition}

\theoremstyle{remark}
\newtheorem{remark}[theorem]{Remark}

\newcommand{\keywords}[1]{\par\medskip\noindent\textbf{Keywords.} #1\par}
\newcommand{\msc}[1]{\par\noindent\textbf{MSC 2020.} #1\par}

\usepackage[colorlinks=true,linkcolor=blue,citecolor=blue,urlcolor=blue,
pdftitle={Hall's condition between adjacent levels of the face lattice of the permutohedron},pdfauthor={Giovanni Rossi}]{hyperref}

\usepackage[capitalize,noabbrev]{cleveref}

\newcommand{\OP}{\mathrm{OP}_n}

\title{Hall's condition between adjacent levels of the face lattice of the permutohedron}
\author{\Large{Giovanni Rossi}\\ Independent researcher, Italy\\ \small{\textit{Email:} \texttt{roxyjean@gmail.com} - ORCID \textit{identifier:} 0000-0001-5215-4325}}
\date{\today}

\begin{document}

\maketitle
\begin{abstract}
Let $\OP$ be the poset of ordered set partitions of $[n]$, ordered by merging adjacent blocks; with a bottom element adjoined it is the face lattice of the permutohedron.
Whether it has the Sperner property is open.
We show that for every pair of adjacent levels, the smaller one can be matched into the larger one by parent--child pairs, except possibly for the pairs $(j,j+1)$ with $(28n-36)/39<j<(13n+5)/18$ (for small $n$ the window is slightly smaller still). This window has width about $n/234\approx0.0043\,n$. The largest level lies near $n/(2\ln2)\approx0.7213\,n$, which we prove using the law of the block sizes of a uniformly random ordered set partition of $[n]$, and for $n$ large it falls inside the window. Since the level sizes are log-concave, Hall's condition for all pairs would make $\OP$ Sperner. The proofs use explicit weights on the parent--child pairs; they reduce to a finite inequality between pairs of adjacent blocks. We explain why the simplest weights fail exactly at the boundary of a much larger window, and report computations supporting Hall's condition inside the window: exactly for $n\le14$, and by equalizing flows for every window pair with $n\le25$.
\end{abstract}
\keywords{Sperner property; ordered set partitions; permutohedron; face lattice; Hall's condition; matchings between adjacent levels; Fubini numbers; renewal theorem}
\msc{Primary 06A07; Secondary 05A18, 05C70, 05D05, 52B05, 60K05}

\section{Introduction}

Let $\OP$ be the set of ordered set partitions $a=(B_1\mid\cdots\mid B_j)$ of $[n]=\{1,\dots,n\}$ into $j\ge1$ non-empty blocks, ordered by merging adjacent blocks.
With a bottom element added it is the face lattice of the permutohedron, in which an element with $j$ blocks is a face of dimension $n-j$ \cite{BilleraSarangarajan,SaganSundaram}.
Its levels are indexed by the number $j$ of blocks, and have sizes $N_j=j!\,S(n,j)$, where $S(n,j)$ is a Stirling number of the second kind.

A graded poset is \emph{Sperner} if its largest antichain is as large as its largest level.
For the lattice $\Pi_n$ of set partitions this is false for $n$ large, by Canfield \cite{Canfield} (see also \cite{CanfieldHarper}).
For $\OP$ the question is open, and we are not aware of a treatment of it in the literature.
We study it through matchings between adjacent levels.

Let $G_j$ be the bipartite graph between levels $j$ and $j+1$ of $\OP$, with an edge for each pair $(a,b)$ in which $b$ is obtained from $a$ by splitting a block into two ordered non-empty parts. We say that \emph{Hall's condition holds for $j$} if $G_j$ has a matching that saturates the smaller of the two levels.

\begin{theorem}\label{thm:main}
Let $1\le j\le n-1$.
\begin{enumerate}[label=(\alph*)]
\item If $j\le 2n/3$ or $39j\le28n-36$, then $G_j$ has a matching saturating level $j$.
\item If $11j\ge8n-2$ or $18j\ge13n+5$, then $G_j$ has a matching saturating level $j+1$.
\end{enumerate}
In particular $N_j\le N_{j+1}$ in case (a) and $N_{j+1}\le N_j$ in case (b).
\end{theorem}

Thus Hall's condition holds for every $j$ outside the window
\begin{equation}\label{eq:window}
\max\Big\{\frac{2n}{3},\frac{28n-36}{39}\Big\}<j<\min\Big\{\frac{8n-2}{11},\frac{13n+5}{18}\Big\}.
\end{equation}
For large $n$ the window is $\big(\tfrac{28}{39}n,\tfrac{13}{18}n\big)$ up to bounded terms, of width $n/234+O(1)$. It contains $n/(2\ln2)\approx0.72135\,n$, the asymptotic location of the largest level (Proposition \ref{prop:mode}), since $28/39=0.7179\ldots$ and $13/18=0.7222\ldots$. It is empty for $n=11,12,14,15,18,25$, equal to $\{7\}$ for $n=10$, to $\{71,72\}$ for $n=100$ and to $\{718,\dots,722\}$ for $n=1000$. The inequalities with $2n/3$ and $(8n-2)/11$ are the ones that matter for $n\le18$ and $n\le91$ respectively.

\begin{conjecture}\label{conj:hall}
Hall's condition holds for every $n$ and every $j$.
\end{conjecture}

\begin{theorem}\label{thm:sperner}
If Conjecture~\ref{conj:hall} holds for $n$, then $\OP$ is Sperner, with width $\max_jN_j$.
\end{theorem}

In Section~\ref{sec:mode} we prove that the largest level lies near $n/(2\ln2)$ and, for $n$ large, inside the window. In Section~\ref{sec:sharp} we show that the two simplest rules cannot be improved, and explain what limits the weighted rules used for the new bounds; in Section~\ref{sec:comp} we report computations: Hall's condition holds for all $j$ for $n\le14$ (including the window pairs $n=10$, $j=7$ and $n=13$, $j=9$), and we find strictly positive equalizing flows for every window pair with $n\le25$.

\paragraph{Levels, children and the reduction to Hall's condition:} an element $a$ with block sizes $s_1,\dots,s_j$ has
\[
d(a)=\sum_{p=1}^j(2^{s_p}-2)
\]
\emph{children} (elements of level $j+1$ obtained by splitting a block $B$ into an ordered pair of non-empty parts: there are $2^{|B|}-2$ such pairs), and each element $b$ of level $j+1$ has exactly $j$ \emph{parents}, obtained by merging two adjacent blocks.

\begin{lemma}\label{lem:lc}
For every $n$ the sequence $N_1,\dots,N_n$ is log-concave, hence unimodal.
\end{lemma}

\begin{proof}
The sequence $S(n,j)$, $1\le j\le n$, is log-concave by Lieb \cite{Lieb}, and the sequence $j!$ is log-concave since ${(j!)^2/((j-1)!\,(j+1)!)=j/(j+1)<1}$. A product of log-concave sequences of positive numbers is log-concave.
\end{proof}

The location of the largest level is determined in Section~\ref{sec:mode}.

\begin{proof}[Proof of Theorem~\ref{thm:sperner}]
By Lemma~\ref{lem:lc} there is $k^*$ with $N_1\le\dots\le N_{k^*}\ge\dots\ge N_n$. For $j<k^*$ take a matching $M_j$ of $G_j$ saturating level $j$, and for $j\ge k^*$ take a matching $M_j$ saturating level $j+1$ (these are the smaller levels).
Every element of a level below $k^*$ lies in exactly one edge of the matching towards $k^*$ and in at most one edge of the matching away from it, and symmetrically above $k^*$.
So the union of the matchings is a disjoint union of paths whose edges are cover relations, each of which contains exactly one element of level $k^*$. These are $N_{k^*}$ chains partitioning $\OP$, and every antichain meets each chain at most once.
\end{proof}

\section{Where the largest level is: compositions and uniform distributions}\label{sec:mode}

The window \eqref{eq:window} was placed around $n/(2\ln2)$. In this section we prove that the largest level is indeed there, and that for $n$ large it lies inside the window, so that the pairs around it are not covered by Theorem~\ref{thm:main}. We first fix the terminology on compositions and the uniform measure, which is the basis of the proofs.

\paragraph{Compositions and ordered set partitions.}
A \emph{composition} of $n$ is a finite sequence $c=(c_1,\dots,c_k)$ of positive integers with $c_1+\dots+c_k=n$; the order of the entries matters, so $(1,2)$ and $(2,1)$ are different compositions of $3$, and $k$ is the number of parts. There are $2^{n-1}$ compositions of $n$. An ordered set partition of $[n]$ is a sequence $(B_1,\dots,B_k)$ of non-empty, pairwise disjoint sets whose union is $[n]$. Its block sizes $(|B_1|,\dots,|B_k|)$ form a composition of $n$, and we say that the ordered set partition has block sizes $c$ when this composition is $c$. For a given composition $c$ the number of ordered set partitions with block sizes $c$ is the multinomial coefficient $n!/(c_1!\cdots c_k!)$: choose the $c_1$ elements of $B_1$, then the $c_2$ elements of $B_2$ among the remaining ones, and so on. Summing over all compositions gives the total number of ordered set partitions, the Fubini number
\[
a_n=\sum_{c}\frac{n!}{c_1!\cdots c_k!}=\sum_{j=1}^{n}j!\,S(n,j),
\]
where $S(n,j)$ is the Stirling number of the second kind, the number of partitions of $[n]$ into $j$ non-empty blocks, without order. The values are $a_n=1,3,13,75,541,\dots$ for $n=1,2,3,4,5,\dots$. The second sum groups the ordered set partitions by number of blocks: $N_j=j!\,S(n,j)$ is the number of those with exactly $j$ blocks, the size of level $j$, so that $a_n=\sum_jN_j=|\OP|$.

\paragraph{The uniform measure.}
A probability measure on a finite set is given by its probability mass function, which assigns to each element a number $\ge0$, with the numbers adding up to $1$; the probability of a subset is the sum of the masses of its elements. The uniform measure on a finite set assigns the same mass, one over the number of elements, to each of them. Here the set is $\OP$, with $a_n$ elements, so the uniform measure $\mathbb P_n$ gives mass $1/a_n$ to each ordered set partition, and a uniform element of $\OP$ is an ordered set partition drawn at random so that all $a_n$ of them are equally likely. Probabilities are then obtained by counting: the probability of an event is the number of ordered set partitions in it divided by $a_n$. In particular the number $k$ of blocks satisfies $\mathbb P_n(k=j)=N_j/a_n$, and a composition $c$ of $n$ has probability $\dfrac{n!}{a_n\,c_1!\cdots c_k!}$, since this is the number of ordered set partitions with block sizes $c$ divided by $a_n$. Two other measures should not be confused with $\mathbb P_n$: the uniform measure on the compositions of $n$ gives each composition the mass $2^{1-n}$, whatever the number of ordered set partitions with those block sizes, and the uniform measure on the levels gives each $j$ the mass $1/n$.

\begin{center}
	\begin{tabular}{l|cccc}
		composition $c$ of $3$&$(3)$&$(1,2)$&$(2,1)$&$(1,1,1)$\\\hline
		ordered set partitions with block sizes $c$&$1$&$3$&$3$&$6$\\
		probability under $\mathbb P_3$&$1/13$&$3/13$&$3/13$&$6/13$\\
		uniform measure on compositions&$1/4$&$1/4$&$1/4$&$1/4$
	\end{tabular}
\end{center}
For $n=3$ one has $a_3=13$, and the three levels have sizes $N_1=1$, $N_2=6$, $N_3=6$, so $\mathbb P_3(k=1)=1/13$ while the uniform measure on levels would give $1/3$.

\paragraph{The law $(p_b)$.}
Put $x_0=\ln2$ and $p_b=x_0^b/b!$ for $b\ge1$. These numbers add up to $e^{x_0}-1=1$, so $(p_b)$ is a probability law on the positive integers, with mean $\mu=\sum_bb\,p_b=x_0e^{x_0}=2\ln2$. Since $\prod_ix_0^{c_i}=x_0^n$, the probability of $c$ under $\mathbb P_n$ equals $\dfrac{n!}{a_n\,x_0^n}\prod_ip_{c_i}$, which is the form used in Proposition \ref{prop:law} below.

\begin{proposition}\label{prop:law}
	Let $X_1,X_2,\dots$ be independent random positive integers, each with law $(p_b)$, and let $S_m=X_1+\dots+X_m$ be their partial sums. We write $\mathbb Q$ for the probability on this sequence, to distinguish it from the uniform measure $\mathbb P_n$ on $\OP$. The sums $S_m$ increase strictly, so at most one of them equals $n$; let $H_n$ be the event that one of them does, and on $H_n$ let $K$ be the index with $S_K=n$. Then $(X_1,\dots,X_K)$ is a composition of $n$ on $H_n$. For every composition $c=(c_1,\dots,c_k)$ of $n$, the probability under $\mathbb P_n$ that a uniform element of $\OP$ has block sizes $(c_1,\dots,c_k)$ in this order equals
	\[
	\frac{n!}{a_n\,x_0^{\,n}}\ \mathbb Q\big((X_1,\dots,X_k)=(c_1,\dots,c_k)\big).
	\]
	Consequently $\mathbb Q(H_n)=a_nx_0^n/n!\to1/\mu$, and the sequence of block sizes of a uniform element has the same law as $(X_1,\dots,X_K)$ conditioned on $H_n$.
\end{proposition}

\begin{proof}
	There are $n!/\prod_ic_i!$ ordered set partitions with block sizes $c_1,\dots,c_k$ in this order, so the probability in the statement is $\frac{n!}{a_n}\prod_i\frac1{c_i!}=\frac{n!}{a_nx_0^n}\prod_ip_{c_i}$, because $\prod_ix_0^{c_i}=x_0^n$; and $\prod_ip_{c_i}=\mathbb Q((X_1,\dots,X_k)=c)$ by independence. The factor in front does not depend on $c$. The event $H_n$ is the disjoint union, over the compositions $c$ of $n$, of the events $\{(X_1,\dots,X_k)=c\}$, so summing over all $c$ gives $1=\frac{n!}{a_nx_0^n}\,\mathbb Q(H_n)$. Since $\sum_b p_b=1$ and $p_1=\ln 2>0$, the law $(p_b)$ is aperiodic, and its mean is $\mu=2\ln 2$.
	Hence, by the renewal theorem \cite[Theorem~1(b), p.~102]{feller1949}, which in this form goes back to Erd\H{o}s, Feller and Pollard \cite{ErdosFellerPollard}, $\mathbb Q(H_n)\to 1/\mu$. Finally $\mathbb Q((X_1,\dots,X_K)=c\mid H_n)=\mathbb Q((X_1,\dots,X_k)=c)/\mathbb Q(H_n)$, which is the $\mathbb P_n$-probability of block sizes $c$.
\end{proof}

\begin{remark}
	The limit $\mathbb Q(H_n)\to1/\mu$ can also be obtained without the renewal theorem. The exponential generating function of the Fubini numbers is $1/(2-e^z)$, that is, the sequence construction applied to $G(z)=e^z-1$. Since $G(x_0)=1$ and $G'(x_0)=2$, the theorem on supercritical sequences \cite[Thm.~V.1, p.~294]{FlajoletSedgewick} (see also \cite[Ch.~6]{Mezo}) gives $a_n/n!=\frac{1}{2x_0^{\,n+1}}\big(1+O(A^n)\big)$ for some $A<1$. Hence $\mathbb Q(H_n)=a_nx_0^n/n!\to1/(2\ln2)=1/\mu$, and the convergence is geometrically fast. For $n=10$ the quantities $a_n/n!$ and $1/(2x_0^{\,n+1})$ already agree to ten decimal places.
\end{remark}

\begin{proposition}\label{prop:mode}
	Let $k$ be the number of blocks of a uniform element of $\OP$. For every $\varepsilon>0$, $\mathbb P_n\big(|k/n-1/(2\ln2)|>\varepsilon\big)\to0$. Moreover, if $j^*=j^*(n)$ is an index with $N_{j^*}=\max_jN_j$, then $j^*/n\to1/(2\ln2)$.
\end{proposition}

\begin{proof}
	Recall that $\mu=2\ln2$. Let $K_n=\#\{m\ge1:S_m\le n\}$. By the strong law for renewal counts \cite[Thm.~2.6.1, p.~79]{durrett2019}, $K_n/n\to1/\mu$ $\mathbb Q$-almost surely. On $H_n$ the number of blocks equals $K_n$, so by Proposition~\ref{prop:law} the law of $k$ under $\mathbb P_n$ is the law of $K_n$ under $\mathbb Q$ conditioned on $H_n$, and
	\[
	\mathbb P_n\big(|k/n-1/\mu|>\varepsilon\big)=\frac{\mathbb Q\big(|K_n/n-1/\mu|>\varepsilon,\ H_n\big)}{\mathbb Q(H_n)}\le\frac{\mathbb Q(|K_n/n-1/\mu|>\varepsilon)}{\mathbb Q(H_n)}\to0,
	\]
	because the numerator tends to $0$ and the denominator to $1/\mu>0$.
	
	For the second claim we may assume $\varepsilon<1/4$, so that all intervals lie in $[0,n]$. Put $c=1/\mu$ and $\pi_j=N_j/a_n=\mathbb P_n(k=j)$. By Lemma~\ref{lem:lc}, $\pi_j$ is non-decreasing for $j\le j^*$ and non-increasing for $j\ge j^*$. Suppose that $j^*\le(c-2\varepsilon)n$ for infinitely many $n$. For such $n$ let $I_1$ be the set of integers in $[(c-2\varepsilon)n,(c-\varepsilon)n]$ and $I_2$ the set of integers in $((c-\varepsilon)n,(c+\varepsilon)n]$. Every element of $I_1$ is at least $j^*$ and is smaller than every element of $I_2$, so $\pi_i\ge\pi_{i'}$ for $i\in I_1$, $i'\in I_2$. Since $|I_1|=\varepsilon n+O(1)$ and $|I_2|=2\varepsilon n+O(1)$, we get $\mathbb P_n(k\in I_1)\ge\big(\tfrac12-o(1)\big)\mathbb P_n(k\in I_2)$. But $\{|k/n-c|\le\varepsilon/2\}\subseteq\{k\in I_2\}$ for $n$ large, so $\mathbb P_n(k\in I_2)\to1$ by the first claim applied with $\varepsilon/2$, while $\mathbb P_n(k\in I_1)\le\mathbb P_n(|k/n-c|\ge\varepsilon)\to0$ by the same claim, a contradiction. The case $j^*\ge(c+2\varepsilon)n$ is symmetric, using the non-decreasing side. Hence $|j^*/n-c|<2\varepsilon$ for $n$ large, for every $\varepsilon>0$.
\end{proof}

Alternatively, by \cite[Thm.~V.1, p.~294]{FlajoletSedgewick} the number of blocks has mean $n/(2\ln2)+O(1)$ and variance $O(n)$ under $\mathbb P_n$, and the first claim follows from Chebyshev's inequality.

\begin{corollary}\label{cor:inside}
	There is $n_0\ge91$ such that for $n\ge n_0$ the index $j^*$ of a largest level satisfies $\frac{28n-36}{39}<j^*<\frac{13n+5}{18}$; that is, the largest level lies inside the window \eqref{eq:window}, which for $n\ge91$ is exactly this interval.
\end{corollary}

\begin{proof}
	One has $28/39=0.71794\ldots<1/(2\ln2)=0.72134\ldots<13/18=0.72222\ldots$, with distances $0.0034$ and $0.00087$ from the endpoints. By Proposition~\ref{prop:mode}, $|j^*/n-1/(2\ln2)|<0.0008$ for $n$ large, and then $\frac{28}{39}n<j^*<\frac{13}{18}n$, which gives both inequalities.
\end{proof}

For moderate $n$ a largest level may sit at the edge of the window or just outside it. Exact values are
\[
\begin{array}{c|cccccc}
	n&20&50&100&200&399&1000\\\hline
	j^*&15&36&72&145&288&722\\
	n/(2\ln2)&14.43&36.07&72.13&144.27&287.82&721.35\\
	\text{window}&(13.4,14.4)&(35.0,36.2)&(70.9,72.5)&(142.7,144.7)&(285.5,288.4)&(717.0,722.5)
\end{array}
\]
so that for $n=20$ and $n=200$ the largest level lies just above the right end of the window, while the level $j^*-1$ lies inside it. An exact computation of all the level sizes for $n\le700$ shows that the last $n$ for which a largest level is not strictly inside the window is $n=613$, where $j^*=443$ is exactly the right endpoint $(13n+5)/18$; for $614\le n\le700$ it is inside, and a floating-point computation gives the same up to $n=4000$. Corollary~\ref{cor:inside} does not give an explicit $n_0$, but these computations suggest $n_0=614$.

\begin{remark}\label{rem:sizes}
By Proposition~\ref{prop:law} and the renewal theorem, for every fixed $m$ the first $m$ block sizes of a uniform element of $\OP$ converge in law, as $n\to\infty$, to $m$ independent variables with law $(p_b)$: $p_1=0.693$, $p_2=0.240$, $p_3=0.0555$, $p_4=0.0096$, $p_5=0.0013$, and blocks of size $6$ or more have total probability $1.7\cdot10^{-4}$; the mean size is $2\ln2=1.386$. So the weights of Section~\ref{sec:outside}, which are tables for blocks of size at most $4$ or $5$ plus a constant for larger blocks, are adapted to the typical element of the middle levels. This does not help inside the window, where Hall's condition has to hold for every composition, not only for typical ones.
\end{remark}

\section{Matchings outside the window: proof of Theorem 1.1}\label{sec:outside}

All the matchings below come from the same principle: it is enough to spread the weight $1$ of each element of $X$ over its neighbours so that no element of $Y$ receives more than $1$ in total, as the following lemma shows.

\begin{lemma}\label{lem:flow}
	Let $X,Y$ be the sides of a bipartite graph. If there are weights $f(x,y)\ge0$ on the edges such that $\sum_yf(x,y)=1$ for every $x\in X$ and $\sum_xf(x,y)\le1$ for every $y\in Y$, then there is a matching saturating $X$.
\end{lemma}
\begin{proof}
If $N(S)$ is the set of neighbours of $S\subseteq X$, then $|S|=\sum_{x\in S}\sum_yf(x,y)=\sum_{y\in N(S)}\sum_{x\in S}f(x,y)\le|N(S)|$, so Hall's condition holds and Hall's theorem gives the matching.
\end{proof}

Hereafter, $a$ denotes an element of level $j$ with block sizes $s_1,\dots,s_j$, and $b$ an element of level $j+1$ with block sizes $b_1,\dots,b_{j+1}$; the parent $a_t$ of $b$ is obtained by merging blocks $t$ and $t+1$ of $b$, for $1\le t\le j$.
The children of $a$ are the elements of level $j+1$ having $a$ as a parent. Each is obtained by choosing a block $p$ of $a$ with $s_p\ge2$ and splitting it into an ordered pair of non-empty sets, a first part of $i$ elements and a second part of $s_p-i$ elements, which replace block $p$ in the sequence. For $1\le p\le j$ and $1\le i\le s_p-1$, the \emph{class} $(p,i)$ of $a$ is the set of children obtained from block $p$ with a first part of $i$ elements. It consists of $\binom{s_p}{i}$ elements, one for each choice of the $i$ elements of the first part, and all of them have the block sizes $(s_1,\dots,s_{p-1},i,s_p-i,s_{p+1},\dots,s_j)$. The classes of $a$ are pairwise disjoint and their union is the set of all children of $a$, because $b$ determines $p$ and $i$. Conversely, $b$ has the $j$ parents $a_1,\dots,a_j$, and in the parent $a_t$ it lies in the class $(p,i)=(t,b_t)$: block $t$ of $a_t$ is the union of blocks $t$ and $t+1$ of $b$, and $b$ arises from $a_t$ by splitting that block with a first part of $b_t$ elements, which is block $t$ of $b$.

\subsection*{Proof of Theorem~\ref{thm:main}(a)}

\paragraph{First condition (the degree argument).}
Let $S$ be a set of elements of level $j$ and $\partial S$ the set of their children.
Since $2^s-2\ge2(s-1)$ for every $s\ge1$, with equality for $s=1,2$, every $a$ in level $j$ satisfies
\[
d(a)\ \ge\ 2\sum_p(s_p-1)=2(n-j)\ \ge\ j\qquad\text{when }j\le 2n/3 .
\]
There are at least $|S|\,j$ edges leaving $S$, and at most $j\,|\partial S|$ edges enter $\partial S$, since each child has $j$ parents. So $|\partial S|\ge|S|$, and Hall's theorem gives the matching.

\paragraph{Second condition.}
Put $w(2)=11$, $w(3)=12$, $w(s)=14$ for $s\ge4$, and $\Omega(1)=0$, $\Omega(s)=(s-1)w(s)$ for $s\ge2$.
Let $U(a)=\sum_p\Omega(s_p)$, which is positive because $j<n$.
Let $a$ send the amount $w(s_p)/\big(\binom{s_p}{i}U(a)\big)$ to each child in the class $(p,i)$.
The total is $\sum_p\sum_{i=1}^{s_p-1}w(s_p)/U(a)=1$.
The amount received by $b$ is
\[
L(b)=\sum_{t=1}^{j}\frac{\kappa(b_t,b_{t+1})}{U(a_t)},\qquad\kappa(x,y)=\frac{w(x+y)}{\binom{x+y}{x}} .
\]
We show $L(b)\le1$ when $39j\le28n-36$; Lemma~\ref{lem:flow} then gives the matching.

\emph{(i) $U(a_t)\ge U(b)$.} Indeed $U(a_t)=U(b)+\Omega(b_t+b_{t+1})-\Omega(b_t)-\Omega(b_{t+1})$, and $\Omega(x+y)\ge\Omega(x)+\Omega(y)$ for all $x,y\ge1$: this is clear if $x=1$ or $y=1$ because $\Omega$ is non-decreasing, and if $x,y\ge2$ then $\Omega(x)+\Omega(y)\le14(x-1)+14(y-1)<14(x+y-1)\le\Omega(x+y)$.

\emph{(ii) A pair inequality.} Let $\varepsilon(2)=3$, $\varepsilon(3)=4$ and $\varepsilon(s)=0$ for $s=1$ and $s\ge4$, so that $\Omega(s)=14(s-1)-\varepsilon(s)$ for all $s\ge1$. For all $x,y\ge1$,
\begin{equation}\label{eq:pairU}
	\kappa(x,y)\ \le\ \frac{\Omega(x)+\Omega(y)}2+\frac{39}2-7(x+y)\ =\ \frac{11}2-\frac{\varepsilon(x)+\varepsilon(y)}2 .
\end{equation}
For $(x,y)=(1,1),(1,2),(1,3)$ the values $\kappa=\tfrac{11}2,4,\tfrac72$ equal the right side, and the same holds for the reversed pairs.
For $(2,2),(2,3),(3,3)$ the left sides are $\tfrac73,\tfrac75,\tfrac7{10}$ and the right sides $\tfrac52,2,\tfrac32$.
All other pairs have $x+y\ge5$ and are not both in $\{2,3\}$, so the right side is at least $\tfrac72$, while $\kappa\le14/(x+y)\le\tfrac{14}5$ because $\binom{x+y}{x}\ge x+y$.

\emph{(iii) Summation.} Sum \eqref{eq:pairU} over the $j$ adjacent pairs of $b$. Using $\sum_t(b_t+b_{t+1})=2n-b_1-b_{j+1}$ and $\sum_t\frac12(\Omega(b_t)+\Omega(b_{t+1}))=U(b)-\frac12(\Omega(b_1)+\Omega(b_{j+1}))$,
\[
\sum_{t=1}^j\kappa(b_t,b_{t+1})\ \le\ U(b)+\frac{39j}2-14n+h(b_1)+h(b_{j+1}),\qquad h(s)=7s-\tfrac12\Omega(s)=7+\tfrac12\varepsilon(s)\le9 .
\]
Hence, by (i), $L(b)\le U(b)^{-1}\sum_t\kappa(b_t,b_{t+1})\le1+\dfrac{39j-28n+36}{2U(b)}\le1$ when $39j\le28n-36$; here $U(b)>0$ because $j+1<n$, which follows from $39j\le28n-36$. \qed

\subsection*{Proof of Theorem~\ref{thm:main}(b)}

\paragraph{First condition (the proportional rule).}
For $b$ put $r_t=\binom{b_t+b_{t+1}}{b_t}^{-1}$ and $Z(b)=\sum_{t=1}^jr_t$.
Let each $b$ send the fraction $r_t/Z(b)$ of a unit of weight to $a_t$. We show that every $a$ receives at most $1$.
Put $e=n-j=\sum_p(s_p-1)$.
Each of the $\binom{s_p}{i}$ children in the class $(p,i)$ sends to $a$ the fraction $\binom{s_p}{i}^{-1}/Z(b)$, where the value of $Z(b)$ depends only on the sizes. So the load on $a$ is
\[
\sum_{p=1}^{j}\ \sum_{i=1}^{s_p-1}\frac1{Z(b_{p,i})},
\]
a sum of $e$ terms, where $b_{p,i}$ is a child in the class $(p,i)$.
For positive integers $x,y$,
\begin{equation}\label{eq:pair}
	\tfrac12-\binom{x+y}{x}^{-1}\le\frac{x+y-2}{6},
\end{equation}
with equality for $(1,1)$, $(1,2)$ and $(2,2)$ (the left side is $0,\tfrac16,\tfrac13$), while for $(1,3)$ the left side is $\tfrac14<\tfrac13$, and for $x+y\ge5$ the right side is at least $\tfrac12$.
Each block of $b$ lies in at most two adjacent pairs, and $\sum_t(b_t-1)=n-(j+1)=e-1$, so
\[
Z(b)\ \ge\ \frac j2-\frac13(e-1)=\frac{3j-2e+2}{6}.
\]
If $j\ge(8n-2)/11$, the denominator is positive, because $3j-2e+2=5j-2n+2\ge(18n+12)/11>0$.
The load on $a$ is therefore at most $6e/(3j-2e+2)$, which is $\le1$ if and only if $8e\le3j+2$, that is, $j\ge(8n-2)/11$.
Lemma~\ref{lem:flow} (with the roles of the levels exchanged) gives a matching saturating level $j+1$.

\paragraph{Second condition.}
Let $u(x,y)=u(y,x)$ be given by
\[
u(1,1)=20,\ u(1,2)=21,\ u(1,3)=20,\ u(1,4)=22,\ u(2,2)=24,\ u(2,3)=20,
\]
and $u(x,y)=0$ if $x+y\ge6$. Put $K(x,y)=u(x,y)/\binom{x+y}{x}$ and $\bar u(s)=\sum_{i=1}^{s-1}u(i,s-i)$, so that
\[
\begin{array}{c|cccccc}
	(x,y)&(1,1)&(1,2)&(1,3)&(1,4)&(2,2)&(2,3)\\\hline
	K&10&7&5&22/5&4&2
\end{array}
\qquad
\begin{array}{c|cccccc}
	s&1&2&3&4&5&\ge6\\\hline
	\bar u(s)&0&20&42&64&84&0
\end{array}
\]
For $b$ put $Y(b)=\sum_{t=1}^jK(b_t,b_{t+1})$.
If $18j\ge13n+5$, then some adjacent pair of $b$ has $b_t+b_{t+1}\le5$, since otherwise $6j\le\sum_t(b_t+b_{t+1})\le2n-2$; so $Y(b)>0$.
Let $b$ send $K(b_t,b_{t+1})/Y(b)$ to $a_t$, a total of $1$.
The $\binom{s_p}{i}$ children of $a$ in the class $(p,i)$ send in total $u(i,s_p-i)/Y(b_{p,i})$, so the load on $a$ is
\[
L(a)=\sum_{p=1}^j\ \sum_{i=1}^{s_p-1}\frac{u(i,s_p-i)}{Y(b_{p,i})} .
\]
We show $L(a)\le1$ when $18j\ge13n+5$.

\emph{(i) $Y(b_{p,i})\ge Y(a)$,} where $Y(a)=\sum_{t=1}^{j-1}K(s_t,s_{t+1})$. Let $c=s_{p-1}$ and $d=s_{p+1}$ (omit the terms that do not exist) and $s=s_p$. Then
\[
Y(b_{p,i})-Y(a)=K(i,s-i)+\big(K(c,i)-K(c,s)\big)+\big(K(s-i,d)-K(s,d)\big)\ \ge0,
\]
because for every fixed $c$ the sequence $K(c,1),K(c,2),\dots$ is non-increasing: it is $10,7,5,\frac{22}5,0,\dots$ for $c=1$; $7,4,2,0,\dots$ for $c=2$; $5,2,0,\dots$ for $c=3$; $\frac{22}5,0,\dots$ for $c=4$; and $0$ for $c\ge5$.

\emph{(ii) A pair inequality.} For all $x,y\ge1$,
\begin{equation}\label{eq:pairD}
	K(x,y)\ \ge\ \frac{\bar u(x)+\bar u(y)}2+36-13(x+y).
\end{equation}
For the six pairs in the table the right sides are $10,7,5,3,4,2$, so there is equality except for $(1,4)$. For $x+y=6$, namely $(1,5),(2,4),(3,3)$ and their reversals, the right side is $0=K$. For $x+y\ge7$ it is $\le0$, because $\bar u(s)\le22(s-1)$ for all $s$ gives $\frac12(\bar u(x)+\bar u(y))\le11(x+y-2)$, and $11(x+y-2)+36-13(x+y)=14-2(x+y)$. The remaining pairs with $x+y\le5$ are the reversals of the six, and $K$ is symmetric.

\emph{(iii) Summation.} Sum \eqref{eq:pairD} over the $j-1$ adjacent pairs of $a$:
\[
Y(a)\ \ge\ \sum_p\bar u(s_p)+36(j-1)-26n+g(s_1)+g(s_j),\qquad g(s)=13s-\tfrac12\bar u(s)\ge13 .
\]
So $Y(a)\ge\sum_p\bar u(s_p)+36j-26n-10\ge\sum_p\bar u(s_p)$ when $18j\ge13n+5$.
By (i), $L(a)\le\sum_p\bar u(s_p)/Y(a)\le1$ (and $L(a)=0$ if all $\bar u(s_p)=0$).
Lemma~\ref{lem:flow} gives a matching saturating level $j+1$. \qed

\section{Sharpness of the rules and the role of the window}\label{sec:sharp}

\begin{proposition}\label{prop:sharp}
\begin{enumerate}[label=(\alph*)]
\item If $n/2\le j\le n-1$, the minimum of $d(a)$ over level $j$ is exactly $2(n-j)$, attained exactly by the elements all of whose blocks have size at most $2$, that is, with $n-j$ blocks of size $2$ and $2j-n$ blocks of size $1$.
\item If $(n+2)/2\le j\le n-1$, the maximum load of the proportional rule in the proof of Theorem~\ref{thm:main}(b) is exactly $6e/(3j-2e+2)$, $e=n-j$. It is attained by every $a$ whose blocks have size at most $2$ and whose first and last blocks have size $1$.
\end{enumerate}
\end{proposition}

\begin{proof}
(a) follows from the inequality $2^s-2\ge2(s-1)$ used in the first proof of Theorem~\ref{thm:main}(a), which is an equality exactly for $s\le2$.
(b) Let $a$ be as stated. A child $b$ of $a$ obtained by splitting a block of size $2$ has all blocks of size at most $2$, and its first and last blocks have size $1$.
Every adjacent pair of $b$ is of type $(1,1)$, $(1,2)$, $(2,1)$ or $(2,2)$, for which \eqref{eq:pair} is an equality, so
$Z(b)=\frac j2-\frac16\sum_t(b_t+b_{t+1}-2)$. Since each block lies in two adjacent pairs, except the first and the last, which have size $1$ and contribute $0$, the sum is $2\sum_t(b_t-1)=2(e-1)$. Hence $Z(b)=(3j-2e+2)/6$ for the children in each of the $e$ classes, and the load on $a$ is $6e/(3j-2e+2)$.
Such an $a$ exists when $e\le j-2$, that is, $j\ge(n+2)/2$, and the upper bound is in the proof of Theorem~\ref{thm:main}(b).
\end{proof}

We checked (b) by enumerating all compositions for six pairs $(n,j)$ with $12\le n\le24$, among them $(16,11)$, $(20,14)$ and $(24,17)$: the maximal load equals the bound in every case, for instance $1.2$ at $(16,11)$ and $1.0769$ at $(24,17)$.
Thus the degree argument and the proportional rule break down exactly at $2n/3$ and $8n/11$: the minimum of $d(a)$ is smaller than $j$ if and only if $j>2n/3$ (for $j\ge n/2$), and the load of the proportional rule exceeds $1$ if and only if $11j<8n-2$ (for $j\ge(n+2)/2$). This shows the limits of these rules, not of Hall's condition itself.
The weighted rules of Theorem~\ref{thm:main} move these boundaries to $\frac{28}{39}n$ and $\frac{13}{18}n$.
The weights were found by linear programming: for a given slope, the pair inequalities \eqref{eq:pairU} and \eqref{eq:pairD} are linear in the weights and in the constants that multiply $x+y$ and the number of blocks, so one can maximize the margin by which they hold.

\paragraph{What the weighted rules can reach.}
The constants $36$ and $5$ in Theorem~\ref{thm:main} are not optimal, but the slopes are the limit of the rules: the exact maxima of the loads over all compositions give, for the rule in (a), the values $0.988$ at $(n,j)=(20,14)$, $1.321$ at $(20,15)$, $0.924$ at $(26,18)$, $1.155$ at $(26,19)$ (here $\frac{28}{39}n=14.4$ and $18.7$), and for the rule in (b) the values $0.940$ at $(22,16)$ and $1.234$ at $(22,15)$ (here $\frac{13}{18}n=15.9$).
Larger supports give better slopes. Let the weights depend on the sizes of the two parts of a block, for blocks of size at most $S$ (and, in case (a), a constant weight for larger blocks). The best slopes of the corresponding linear programs, computed in floating point and not certified, are as follows (in case (b) with $S\ge8$ the condition (i) of the proof is imposed in its exact form, not through monotonicity of $K(c,\cdot)$):
\[
\begin{array}{c|cccc|cccc}
&\multicolumn{4}{c|}{\text{lower slope, case (a)}}&\multicolumn{4}{c}{\text{upper slope, case (b)}}\\
S&3&4&6&8&5&6&7&8\\\hline
\text{slope}&0.71795&0.72022&0.72097&0.72115&0.72222&0.72149&0.72137&0.72135
\end{array}
\]
The mode is at $1/(2\ln2)=0.721348\ldots$. The values stabilize: in our linear programs the slope in case~(a) equals $75/104=0.72115\ldots$ for every $S\ge8$ (checked up to $S=16$), and in case~(b) it equals $0.72135\ldots$ for every $S\ge10$ (checked up to $S=14$). Hence rules of this kind reach Hall's condition only outside a window of width about $2\cdot10^{-4}n$ around the mode $n/(2\ln2)$, and closing the gap would need rules of a different kind.

Inside a window of width $o(n)$ such rules are not expected to work. Put $\rho_j=\min(N_j,N_{j+1})/\max(N_j,N_{j+1})$, the mean load of a flow that saturates the smaller level. Then $\rho_j$ is very close to $1$ near the largest level; for instance $\rho_{15}=0.9902$ for $n=21$. A proof of Conjecture~\ref{conj:hall} there must produce almost perfectly equalizing flows, and in the linear programs above we observed that the worst load of a rule with weights of bounded support exceeds its mean at the mode.

\section{Computations}\label{sec:comp}

The graph computations are exact integer computations; the equalizing flows are in floating point.

\paragraph{The graphs for small $n$.}
For $n\le8$ we built the graphs $G_j$ in full and computed maximum matchings: Hall's condition holds for every $j$.

\paragraph{The quotient by $S_n$.}
The symmetric group acts on $\OP$ by permuting labels, and the orbits are indexed by the compositions of $n$, the orbit of a composition $c=(c_1,\dots,c_j)$ having size $w(c)=n!/\prod_ic_i!$.
Since $S_n$ is transitive on each orbit, the edges between two orbits of adjacent levels form a biregular bipartite graph. Hence $G_j$ has a matching saturating the smaller level if and only if the flow problem on the quotient has a solution, in which each orbit of the smaller level supplies $w(c)$, each orbit of the larger level has capacity $w(c')$, and arcs join compositions related by splitting a part. (A flow on the quotient lifts to a fractional flow on $G_j$ by spreading it uniformly along the edges between two orbits, and a fractional flow of the right value gives a matching by Lemma~\ref{lem:flow}; conversely, averaging a matching over $S_n$ gives a flow on the quotient.)
Exact integer maximum flows on the quotient show that Hall's condition holds for every $j$, for $3\le n\le14$. The window \eqref{eq:window} is non-empty for $n=10$ ($j=7$) and $n=13$ ($j=9$), so Conjecture~\ref{conj:hall} holds in the window for these $n$.
The \emph{normalized} matching property, in which each element of the smaller level $S$ sends weight $1/|S|$ and each element of the larger level $T$ receives at most $1/|T|$, holds for every $j\ge3$ in this range, and fails for $j=2$ for $n=12,13,14$, where Hall's condition nevertheless holds.

\paragraph{Equalizing flows in the window.}
Let $\rho=N_{\rm small}/N_{\rm large}$ for a pair of adjacent levels. A fractional flow is \emph{equalizing} if each element of the smaller level sends weight $1$ and each element of the larger level receives exactly $\rho$.
If a strictly positive equalizing flow exists on the quotient, Hall's condition holds for the pair, because $\rho<1$.
We computed equalizing flows by iterative scaling of the adjacency matrix of the quotient (that is, Sinkhorn's algorithm with marginals $w(\cdot)$ on the smaller level and $\rho\,w(\cdot)$ on the larger one) for every window pair with $n\le25$, namely $(n,j)=(10,7),(13,9),(16,11),(17,12),(19,13),(20,14),(21,15),(22,15),(23,16),(24,17)$ (the window is empty for the other $n\le25$ with $n\ge10$).
In all cases the iteration converged, with maximal relative deviation of the row sums below $10^{-10}$ and all scaling factors positive and finite.
Values of $\rho$: $0.980$ for $(n,j)=(10,7)$, $0.875$ for $(16,11)$, $0.929$ for $(20,14)$, $0.990$ for $(21,15)$, $0.968$ for $(24,17)$.
These are floating-point computations and are not certified.

\section{Open questions}

\begin{enumerate}[label=\arabic*.]
\item \textbf{The window.} Prove Conjecture~\ref{conj:hall} for $j$ in the window \eqref{eq:window}. The weight tables of Section~\ref{sec:sharp} stop improving at slopes $0.72115$ and $0.72135$, so they leave a window of width about $2\cdot10^{-4}n$ around the mode; closing it needs rules of a different kind, and inside a window of width $o(n)$ explicit, almost perfectly equalizing flows.
\item \textbf{The normalized property.} For $3\le n\le18$ the normalized matching property holds for all $j\ge3$, but it fails for $j=2$ from $n=12$ on (checked up to $n=45$), for $j=3$ from $n=20$ (up to $36$) and for $j=4$ from $n=28$ (up to $30$). For which pairs $(n,j)$ does it hold? In particular, does it hold for all $j$ in a neighbourhood of the mode?
\item \textbf{Sperner property.} Is $\OP$ Sperner for every $n$? By Theorem~\ref{thm:sperner} this follows from Conjecture~\ref{conj:hall}. The contrast with $\Pi_n$, which is not Sperner for large $n$ \cite{Canfield,CanfieldHarper}, makes the question natural.
\end{enumerate}

\section*{Acknowledgment}
The author gratefully acknowledges the assistance of Claude (Anthropic) in computations and verification, while retaining full responsibility for the whole manuscript.

\end{document}